\documentclass[12pt]{article}
\usepackage{hyperref}
\usepackage{amsfonts}
\usepackage{enumerate}
\begin{document}

\begin{center}
{\large\bf Solutions of Higher Order Homogeneous Linear Matrix Differential Equations: Singular Case}

\vskip.20in

Grigoris Kalogeropoulos$^{1}$, Charalambos P. Kontzalis$^{2}$\\[2mm]
{\footnotesize
$^{1}$Department of Mathematics, University of Athens, Greece\\[5pt]
$^{2}$Department of Informatics, Ionian University, Corfu, Greece}
\end{center}

{\footnotesize
\noindent
\textbf{Abstract:}  The main objective of this talk is to develop a matrix pencil approach for the study of an initial value problem of a class of singular linear matrix differential equations whose coefficients are constant matrices. By using matrix pencil theory we study the cases of non square matrices and of square matrices with an identically zero matrix pencil. Furthermore we will give necessary and sufficient conditions for existence and uniqueness of
solutions and we will see when the uniqueness of solutions is not valid.  Moreover we
provide a numerical example.\\
\\[3pt]

\vskip.2in

\section{Introduction}

Linear Matrix Differential Equations (LMDEs) are
inherent in many physical, engineering, mechanical, and
financial/actuarial models. In finance for instance, we provide the well-known
input-output Leondief model and its several important extensions, see [1]. Describing electrical circuits and impulsive behavior is another application described in [1]. A second order LMDE which describes the vibration of a building is another example, see [43]. For other applications of continuous-discrete time systems see [32-35, 44-48]. We consider 
\begin{equation}
A_nX^{(n)}(t)+A_{n-1}X^{(n-1)}(t)+...+A_1X'(t)+A_0X(t)=0_{m, 1} 
\end{equation}
with known initial conditions
\begin{equation}
X(t_0),X'(t_0),...,X^{n-1}(t_0)
\end{equation}
where  $A_i, i=0,1,...,n \in \mathcal{M}({m \times r;\mathcal{F}})$, (i.e. the
algebra of square matrices with elements in the field
$\mathcal{F}$) with  $X(t) \in\mathcal{M}({m \times
1;\mathcal{F}})$ and det$A_n$=0 if $A_n$ is square. For the sake of simplicity we set
${\mathcal{M}}_m  = {\mathcal{M}}({m \times m;\mathcal{F}})$ and
${\mathcal{M}}_{mr}  = {\mathcal{M}}({m \times r;\mathcal{F}} )$. In the sequel we adopt the following notations 
\begin{equation}
    \begin{array}{c}
    Y_1 (t)=X (t),\\
    Y_2 (t)=X' (t),\\
    \dots\\
    Y_{n-1}(t)=X^{(n-1)}(t),\\
    Y_n (t)=X^{(n-1)} (t).
    \end{array}
   \end{equation}
   \begin{equation}
   \begin{array}{c}
    Y_1' (t)=X' (t)=Y_2(t),\\
    Y_2' (t)=X'' (t)=Y_3(t),\\
    \dots\\
    Y_{n-1}(t)=X^{(n-1)}(t)=Y_n(t),\\
    A_nY_n' (t)=A_nX^{(n)}(t)=-A_{n-1}Y_n(t)-...-A_1Y_2(t)-A_0Y_1(t).
    \end{array}
    \end{equation}
Or in Matrix form
\begin{equation}
FY'(t)=GY(t)
\end{equation}
with known initial conditions
\begin{equation}
Y(t_0)
\end{equation}        
where $Y(t) =[ Y_1^T (t)Y_2^T(t)\dots Y_n^T
(t)]^T$ (where $()^T$ is the transpose tensor) and the
coefficient matrices $F, G$ are given by
\begin{equation}
 F = \left[
\begin{array}{ccccc} I_m&0_{m, m}&...&0_{m, m}&0_{m, m}\\0_{m, m}&I_m&...&0_{m, m}&0_{m, m}\\\vdots&\vdots&\ddots&\vdots&\vdots\\0_{m, m}&0_{m, m}&...&I_m&0_{m, m}\\0_{m, r}&0_{m, r}&...&0_{m, r}&A_n
\end{array}
\right] 
 G = \left[\begin{array}{cccc} 0_{m, m}&I_m&\ldots&0_{m, m}\\0_{m, m}&0_{m, m}&\ldots&0_{m, m}\\\vdots&\vdots&\ddots&\vdots\\0_{m, m}&0_{m, m}&\ldots&I_m\\-A_0&-A_1&\ldots&-A_{n-1}\end{array}\right]
\end{equation}
with corresponding dimension of F, G and $Y(t)$
$mn\times (mn+r-m)$, $mn\times (mn+r-m)$ and $(mn+r-m)\times 1$, respectively. The matrix F is singular if $A_n$ (and consequently F) is singular.

The matrix pencil theory has been extensively used for the study of
Linear Matrix Differential/Difference Equations with time invariant
coefficients, see for instance [1-6, 8, 10, 11, 13-16, 19-23, 27-29, 36-42]. Recently this has been extended also to Linear Matrix Fractional Differential/Difference Equations, see [21, 22, 30, 31].
A matrix pencil is a family of matrices sF-G, parametrized by a complex number s, see [5, 6, 9, 20]. When G is square and $F=I_n$, where $I_n$ is the identity matrix, the zeros of the function det(sF-G) are the eigenvalues of G. Consequently, the problem of finding the nontrivial solutions of the equation
\begin{equation}
sFX=GX
\end{equation}
is called the generalized eigenvalue problem. Although the generalized eigenvalue problem looks like a simple generalization of the usual eigenvalue problem, it exhibits some important differences.

In the first place, if F, G are square matrices (for r=m) it is possible for F to be singular, in which case the problem has infinite eigenvalues. To see this, write the generalized eigenvalue problem in the reciprocal form
\begin{equation}
FX=s^{-1}GX
\end{equation}
If F is singular with a null vector X, then $GX=0_{mn, 1}$, so that X is an eigenvector of the reciprocal problem corresponding to eigenvalue $s^{-1}=0$; i.e., $s=\infty$. In this case the solutions of system (1) have been fully analyzed in [1-3, 5, 6, 14, 19, 20, 22, 23].

Second, it is possible for det(sF-G) to be identically zero, independent of s. Finally it is possible both matrices F, G to be non square (for r$\neq$m). In this article we will consider these last two cases.

\section{Singular matrix pencils: Mathematical background and notation}
In this section we will give the mathematical background and the notation that is used throughout the paper
\\\\
\textbf{Definition 2.1}
Given $F,G\in \mathcal{M}_{nm}$  and an indeterminate $s\in\mathcal{F}$, the matrix
pencil $sF-G$ is called regular when  $m=n$ and  $\det(sF-G)\neq 0$.
In any other case, the pencil will be called singular.
\\\\
\textbf{Definition 2.2}
The pencil $sF-G$  is said to be \emph{strictly equivalent} to the
pencil $s\tilde F - \tilde G$ if and only if there exist nonsingular
$P\in\mathcal{M}_n$ and $Q\in\mathcal{M}_m$ such as
\[
    P({sF - G})Q = s\tilde F - \tilde G.
\]
In this article, we consider the case that the pencil is \emph{singular}. The main results for this case
are analytically presented. Unlike the case of the regular pencils, the characterisation of a singular matrix pencil, apart from the set of the determinantal divisors requires the definition of additional sets of invariants, the minimal indices. The distinguishing feature of a singular pencil sF-G is that either $mn\neq mn$ or mn=mn and sF-G$\equiv 0$. Let $\mathcal{N}_r$, $\mathcal{N}_l$ be right, left null space of a matrix respectively. Then the equations
\begin{equation}
(sF-G)U(s)=0_{m, 1}
\end{equation}
\begin{equation}
V^T(s)(sF-G)=0_{1, m}
\end{equation}
have solutions in U(s), V(s), which are vectors in the rational vector spaces $\mathcal{N}_r(sF-G)$ and $\mathcal{N}_l(sF-G)$ respectively. The binary vectors X(s) and $Y^T(s)$ express dependence relationships among the colums or rows of sF-G respectively. U(s), V(s) is polynomial vectors. Let d=dim$\mathcal{N}_r(sF-G)$ and t=$\mathcal{N}_l(sF-G)$. It is known [5, 6, 9] that $\mathcal{N}_r(sF-G)$, $\mathcal{N}_l(sF-G)$, as rational vector spaces, are spanned by minimal polynomial bases of minimal degrees 
\begin{equation}
\epsilon_1=\epsilon_2=...=\epsilon_g=0<\epsilon_{g+1}\leq...\leq\epsilon_d
\end{equation}
and 
\begin{equation}
\zeta_1=\zeta_2=...=\zeta_h=0<\zeta_{h+1}\leq...\leq\zeta_t
\end{equation}
respectively. The set of minimal indices $\epsilon_i$ and $\zeta_j$ are known [5, 6, 9] as \textit{column minimal indices} (c.m.i.) and \textit{row minimal indices} (r.m.i) of sF-G respectively. 
To sum up in the case of a singular matrix pencil, we have invariants, a set of \emph{elementary divisors} (e.d.) and \emph{minimal indices}, of the following type:

\begin{itemize}
    \item e.d. of the type  $(s-a)^{p_j}$, \emph{finite elementary
    divisors} (nz. f.e.d.)
    \item e.d. of the type  $\hat{s}^q=\frac{1}{s^q}$, \emph{infinite elementary divisors}
    (i.e.d.).
    \item m.c.i. of the type $\epsilon_1=\epsilon_2=...=\epsilon_g=0<\epsilon_{g+1}\leq...\leq\epsilon_d$, \emph{minimal column indices}
    \item m.r.i. of the type $\zeta_1=\zeta_2=...=\zeta_h=0<\zeta_{h+1}\leq...\leq\zeta_t$, \emph{minimal row indices}
\end{itemize}
\textbf{Definition 2.3.} Let $B_1 ,B_2 ,\dots, B_n $ be elements of $\mathcal{M}_n$. The direct sum
of them denoted by $B_1  \oplus B_2  \oplus \dots \oplus B_n$ is
the blockdiag$\left[\begin{array}{cccc} B_1& B_2& \dots& B_n\end{array}\right]$.
\\\\
The existence of a complete set of invariants for singular pencils implies the existence of canonical form, known as Kronecker canonical form [5, 6, 9] defined by 
\begin{equation}
sF_K  - Q_K :=
sI_p  - J_p  \oplus sH_q  - I_q \oplus sF_{\epsilon}-G_{\epsilon}\oplus sF_{\zeta}-G_{\zeta}\oplus 0_{h, g}
\end{equation}
where $sI_p  - J_p$, $sH_q  - I_q$ are defined in section 2. Let $B_1 ,B_2 ,\dots, B_n $ be elements of $\mathcal{M}_n$. The direct sum
of them denoted by $B_1  \oplus B_2  \oplus \dots \oplus B_n$ is the {block diag}$\{ B_1 ,B_2 ,\dots, B_n \}$.

where $sI_p  - J_p$ is uniquely defined by the set of f.e.d.
\begin{equation}
  ({s - a_1 })^{p_1 } , \dots ,({s - a_\nu  }
 )^{p_\nu },\quad \sum_{j = 1}^\nu  {p_j  = p}
\end{equation}
of sF-G  and has the form
\begin{equation}
    sI_p  - J_p  := sI_{p_1 }  - J_{p_1 } (
    {a_1 }) \oplus  \dots  \oplus sI_{p_\nu  }  - J_{p_\nu  }
    ({a_\nu  }) 
\end{equation}
The q  blocks of the second uniquely defined block
$sH_q -I_q$ correspond to the i.e.d.
\begin{equation}
  \hat s^{q_1} , \dots ,\hat s^{q_\sigma}, \quad \sum_{j =
  1}^\sigma  {q_j  = q}
\end{equation}
of sF-G  and has the form
\begin{equation}
    sH_q  - I_q  := sH_{q_1 }  - I_{q_1 }  \oplus
    \dots  \oplus sH_{q_\sigma  }  - I_{q_\sigma}
\end{equation}
Thus, $H_q$  is a nilpotent element of $\mathcal{M}_n$  with index
$\tilde q = \max \{ {q_j :j = 1,2, \ldots ,\sigma } \}$, where
\[
    H^{\tilde q}_q=0_{q, q},
\]
and $I_{p_j } ,J_{p_j } ({a_j }),H_{q_j }$ are defined as
\begin{equation}
   I_{p_j }  = \left[\begin{array}{ccccc} 
   1&0& \ldots & 0&0\\
   0& 1 &  \ldots&0 &0 \\
   \vdots & \vdots & \ddots & \vdots &\vdots \\
   0 & 0 & \ldots  & 0 &1
   \end{array}\right]
   \in {\mathcal{M}}_{p_j } , 
   \end{equation}
   \begin{equation}
   J_{p_j } ({a_j }) =  \left[\begin{array}{ccccc}
   a_j  & 1 & \dots&0  & 0  \\
   0 & a_j  &   \dots&0  & 0  \\
    \vdots  &  \vdots  &  \ddots  &  \vdots  &  \vdots   \\
   0 & 0 &  \ldots& a_j& 1\\
   0 & 0 & \ldots& 0& a_j
   \end{array}\right] \in {\mathcal{M}}_{p_j }
\end{equation}
\begin{equation}
 H_{q_j }  = \left[
\begin{array}{ccccc} 0&1&\ldots&0&0\\0&0&\ldots&0&0\\\vdots&\vdots&\ddots&\vdots&\vdots\\0&0&\ldots&0&1\\0&0&\ldots&0&0
\end{array}
\right] \in {\mathcal{M}}_{q_j }.
  \end{equation}
For algorithms about the computations of the jordan matrices see [5, 6, 8, 9, 11, 18, 24, 25]. For the rest of the diagonal blocks of $F_K$ and $G_K$, s$F_{\epsilon}-G_{\epsilon}$ and s$F_{\zeta}-G_{\zeta}$,  the matrices $F_{\epsilon}$, $G_{\epsilon}$ are defined as
\begin{equation} 
F_\epsilon=blockdiag\left\{L_{\epsilon_{g+1}}, L_{\epsilon_{g+2}}, ..., L_{\epsilon_d}\right\}
\end{equation}
 Where $L_\epsilon= \left[
\begin{array}{ccc} I_\epsilon & \vdots & 0_{\epsilon, 1}
\end{array}
\right]$, for $\epsilon=\epsilon_{g+1}, ..., \epsilon_d$
\begin{equation}
G_\epsilon=blockdiag\left\{\bar L_{\epsilon_{g+1}}, \bar L_{ \epsilon_{g+2}}, ..., \bar L_{\epsilon_d}\right\}
\end{equation}  
Where $\bar L_\epsilon=\left[
\begin{array}{ccc} 0_{\epsilon, 1} & \vdots & I_\epsilon
\end{array}
\right]$, for $\epsilon=\epsilon_{g+1}, ..., \epsilon_d$. The matrices $F_{\zeta}$, $G_{\zeta}$ are defined as
\begin{equation} 
F_\zeta=blockdiag\left\{L_{\zeta_{h+1}}, L_{\zeta_{h+2}}, ..., L_{\zeta_t}\right\}
\end{equation}
 Where $L_\zeta= \left[
\begin{array}{c} I_\zeta \\ 0_{1, \zeta}
\end{array}
\right]$, for $\zeta=\zeta_{h+1}, ..., \zeta_t$
\begin{equation}
G_\zeta=blockdiag\left\{\bar L_{\zeta_{h+1}}, \bar L_{ \zeta_{h+2}}, ..., \bar L_{\zeta_t}\right\}
\end{equation}
 Where $\bar L_\zeta=\left[
\begin{array}{c} 0_{1, \zeta}\\I_\zeta
\end{array}
\right]$, for $\zeta=\zeta_{h+1}, ..., \zeta_t$

\section{Solution space form of LMDEs with singular matrix pencil}

Following the above given analysis, there exist non-singular matrices P, Q such that 
\begin{equation}
PFQ=F_K
\]
\[
PGQ=G_K
\end{equation}
Let 
\begin{equation}
Q=\left[\begin{array}{ccccc}Q_p & Q_q &Q_\epsilon & Q_\zeta & Q_g\end{array}\right]
\end{equation}
where $Q_p\in \mathcal{M}_{(mn)p}$, $Q_q\in \mathcal{M}_{(mn)q}$, $Q_\epsilon\in \mathcal{M}_{(mn)\epsilon}$, $Q_\zeta\in \mathcal{M}_{(mn)\zeta}$ and $Q_g\in \mathcal{M}_{(mn)g}$
\\\\
\textbf{Lemma 3.1.}
System (5) is divided into five subsystems:
\begin{equation}
    Z_p'(t) = J_p Z_p(t) 
\end{equation}
the subsystem
\begin{equation}
    H_q Z_q'(t) = Z_q(t)
\end{equation}
the subsystem
\begin{equation}
    F_\epsilon Z_\epsilon'(t)=G_\epsilon Z_\epsilon(t)
\end{equation}
the subsystem
\begin{equation}
    F_\zeta Z_\zeta'(t)=G_\zeta Z_\zeta(t)
\end{equation}
the subsystem
\begin{equation}
    0_{h, g}\cdot Z_g'(t)=0_{h, g}\cdot Z_g(t)
\end{equation}
\textbf{Proof.}
Consider the transformation
\begin{equation}
    Y(t)=QZ(t)
\end{equation}
Substituting the previous expression into (5) we obtain
\[
    FQZ'(t)=GQZ(t).
\]
Whereby, multiplying by  P and using (14), (26), we arrive at
\begin{equation}
    F_KZ'(t)=G_K Z(t).
\end{equation}
Moreover, we can write Z(t) as
\[
Z(t)=\left[
\begin{array}{c} Z_p(t)\\Z_q(t)\\Z_\epsilon(t)\\Z_\zeta(t)\\Z_g(t)
\end{array}
\right]
\]
where $Z_p(t)\in \mathcal{M}_{p1}$, $Z_q(t) \in \mathcal{M}_{q1}$, $Z_\epsilon(t)\in \mathcal{M}_{p1}$, $Z_\zeta(t) \in \mathcal{M}_{q1}$ and $Z_g(t)\in \mathcal{M}_{p1}$, $Z_q(t) \in \mathcal{M}_{q1}$. Taking into account the above expressions, we arrive easily at
(28-32).
\\\\
Solving the system (5) is equivalent solving subsystems (28-32).
\\\\
\textbf{Remark 3.1.} System (28) is a regular type system and its solution is given from, see [1-3, 7, 8, 10, 19, 20] 
\begin{equation}
Z_p(t)=e^{J_p(t-t_0)}Z_p(t_0)
\end{equation}
\textbf{Remark 3.2.} System (29) is a singular type system but its solution is very ease to compute, see [5, 6, 14-16, 22, 23]
\begin{equation}
Z_q(t)=0_{q,1}
\end{equation}
\textbf{Proposition 3.1.} The subsystem (30) has infinite solutions
\\\\
\textbf{Proof.} If we set 
\[
Z_\epsilon(t)=\left[\begin{array}{c} Z_{\epsilon_{g+1}}(t)\\Z_{\epsilon_{g+2}}(t)\\\vdots\\Z_{\epsilon_d}(t)\end{array}\right]
\]
The system (30) can be written as:
\[
blockdiag\left\{L_{\epsilon_{g+1}}, ..., L_{\epsilon_d}\right\}\left[\begin{array}{c} Z_{\epsilon_{g+1}}'(t)\\Z_{\epsilon_{g+2}}'(t)\\\vdots\\Z_{\epsilon_d}'(t)\end{array}\right]=blockdiag\left\{\bar L_{\epsilon_{g+1}}, ..., \bar L_{\epsilon_d}\right\}\left[\begin{array}{c} Z_{\epsilon_{g+1}}(t)\\Z_{\epsilon_{g+2}}(t)\\\vdots\\Z_{\epsilon_d}(t)\end{array}\right]
\]
Then for the non-zero blocks we have:
\begin{equation}
\begin{array}{ccc} L_\epsilon Z_\epsilon'(t)=\bar L_\epsilon Z_\epsilon(t) & , & \epsilon=\epsilon_{g+1}, ..., \epsilon_d \end{array}
\end{equation}
Using (22), (23) a typical equation from (37) can be written as
\begin{equation}
\left[\begin{array}{ccc} I_\epsilon & \vdots & 0_{\epsilon, 1}\end{array}\right]Z_\epsilon'(t)=\left[
\begin{array}{ccc} 0_{\epsilon, 1} & \vdots & I_\epsilon
\end{array}
\right]Z_\epsilon(t)
\end{equation}
Every solution of the system (37) can be assigned arbitrary. The system (38) is of a regular type differential system. It is clear from the above analysis that in every one of the d-g subsystems one of the coordinates of the solution has to be arbitrary by assigned total. The system has no unique solution and can be taken arbitrary
\begin{equation}
Z_\epsilon(t)=C_{k,1}
\end{equation}
\textbf{Proposition 3.2.} The subsystem (31) has the unique solution
\begin{equation}
Z_\zeta(t)=0_{g, 1}
\end{equation}
\textbf{Proof.} If we set 
\[
Z_\zeta(t)=\left[\begin{array}{c} Z_{\zeta_{h+1}}(t)\\Z_{\zeta_{h+2}}(t)\\\vdots\\Z_{\zeta_t}(t)\end{array}\right]
\]
The system (31) can be written as:
\[
blockdiag\left\{L_{\zeta_{h+1}}, ..., L_{\zeta_t}\right\}\left[\begin{array}{c} Z_{\zeta_{h+1}}'(t)\\Z_{\zeta_{h+2}}'(t)\\\vdots\\Z_{\zeta_t}'(t)\end{array}\right]=blockdiag\left\{\bar L_{\zeta_{h+1}}, ..., \bar L_{\zeta_t}\right\}\left[\begin{array}{c} Z_{\zeta_{h+1}}(t)\\Z_{\zeta_{h+2}}(t)\\\vdots\\Z_{\zeta_t}(t)\end{array}\right]
\]
Then for the non-zero blocks we have:
\begin{equation}
\begin{array}{ccc} L_\zeta Z_\zeta'(t)=\bar L_\zeta Z_\zeta(t) & , & \zeta=\zeta_{h+1}, ..., \zeta_t \end{array}
\end{equation}
Because of structure of $(L_\zeta, \bar L_\zeta)$ blocks, it is readily shown that the only solution of (41) is the zero solution.
\begin{equation}
Z_\zeta(t)=0_{t-h, 1}
\end{equation}
\textbf{Remark 3.3.}The last system (32) has an infinite number of solutions that can be taken arbitrary
\begin{equation}
Z_g(t)=C_{k,2}
\end{equation}
\textbf{Theorem 3.1.} Consider the system (5), with known initial conditions (6) and  let the matrix pencil sF-G be singular. Then the solution is unique if and only if the c.m.i. are zero
\begin{equation}
dim\mathcal{N}_r(sF-G)=0
\end{equation}
and
\begin{equation}
Y(t_0)\in colspan Q_p
\end{equation}
The unique solution is then given from the formula
\begin{equation}
    Y(t)=Q_pe^{J_p(t-t_0)}Z_p(t_0) 
\end{equation}
where $Z_p(t_0)$ is the unique solution of the algebraic system $Y(t_0)=Q_pZ_p(t_0)$. In any other case the system has infinite solutions.
\\\\
\textbf{Proof.} First we consider that the system has non zero c.m.i and non zero r.m.i. Consider the transformation (33), then from lemma 3.1, remarks 3.1, 3.2, 3.3 and propositions 3.1, 3.2, the system (5) is divided into the subsystems (28-32) with solutions (35), (36), (39), (40), (43) respectively, then
\[
     Y(t) = QZ(t) =
     \left[\begin{array}{ccccc}Q_p & Q_q &Q_\epsilon & Q_\zeta & Q_g\end{array}\right]
     \left[\begin{array}{c}
     e^{J_p(t-t_0)}Z_p(t_0)  \\
     0_{q, 1}\\C_{k,1}
     \\0_{t-h, 1}
     \\C_{k,2}
     
    \end{array}\right] 
    \]
    \[
    Y(t) =
    Q_pe^{J_p(t-t_0)}Z_p(t_0)+Q_\epsilon C_{k,1}+Q_g C_{k,2}
\]
Since $C_{k,1}$ and $C_{k,2}$ can be taken arbitrary, it is clear that the general singular LMDE for every suitable defined initial condition has an infinite number of solutions, so it doesn't represent a dynamical system.
 It is clear that the existence of c.m.i. is the reason that ths system (30) and consequently system (32) exist. These systems as shown in propositions 3.1 and remark 3.3 have always infinite solutions. Thus a necessary condition for the system to have unique solution is  not to have any c.m.i. which is equal to 
\[
dim\mathcal{N}_r(sF-G)=0
\]
In this case the Kronecker canonical form of the pencil sF-G has the following form
\begin{equation}
sF_K  - Q_K :=
sI_p  - J_p  \oplus sH_q  - I_q \oplus sF_{\zeta}-G_{\zeta}
\end{equation}   
and then the system (5) is divided into the three subsystems (28), (29), (31) with solutions (35), (36), (40) respectively. Thus
\[
     Y(t) = QZ(t) =
     \left[\begin{array}{ccc}Q_p & Q_q & Q_\zeta\end{array}\right]
     \left[\begin{array}{c}
     e^{J_p(t-t_0)}Z_p(t_0)  \\
     0_{q, 1}
     \\0_{t-h, 1}
     \end{array}\right] 
    \]
    \[
Y(t) =
    Q_pe^{J_p(t-t_0)}Z_p(t_0)
\]
The solution that exists if and only if
\[
Y(t_0)=Q_pZ_p(t_0)
\]
 or
\[
Y(t_0)\in colspan Q_p
\]
In this case the system has the unique solution 
\[
   Y(t)=Q_pe^{J_p(t-t_0)}Z_p(t_0)
\]
\textbf{Remark 2.2.5.}
It follows that the system (1) with known initial conditions (2) has a unique solution if and only if the initial value problem (5), (6) has a unique solution. Then its analytic solution is given by
\begin{equation}
    X(t)=Q_p^1e^{J_pt}Z_p(t_0) 
\end{equation}
where $Q_p^1$ is defined from the matrix 
\[
     Q_p =
     \left[\begin{array}{c}
     {Q_p^1}  \\
     {Q_p^2}
    \end{array}\right]  .
\]
$Q_p^1\in \mathcal{M}_{pp}$
and where $Q^1$ is defined from the matrix 
\[
     Q =
     \left[\begin{array}{c}
     {Q^1}  \\
     {Q^2}
    \end{array}\right]  
\]
and $Q^1\in \mathcal{M}_{(mn)(mn)}$

\section{Numerical example}
\subsection{Example 1}
Consider the matrix differential equation (5) and let 
\[
F=\left[\begin{array}{ccccccc} 2&1&1&0&0&0&0\\1&3&1&1&0&0&0\\1&1&2&1&0&0&0\\0&1&1&1&0&0&0\\0&0&0&0&0&0&0\\0&0&0&0&1&0&0\\0&1&0&0&0&0&1\end{array}\right],
\]
and
\[
G=\left[\begin{array}{ccccccc} 1&1&1&0&0&0&1\\0&3&2&2&0&1&1\\1&2&3&2&0&0&0\\0&2&2&2&0&0&0\\0&0&0&0&1&0&0\\0&0&0&0&0&0&0\\0&0&0&0&0&1&0\end{array}\right]
\]
Then det[sF-G]=0. The invariants of the pencil are,  s-2, s-1 the finite elementary divisors, $\epsilon_1$=0, $\epsilon_2$=2 the c.m.i and 
$\zeta_1$=0, $\zeta_2$=1 are the r.m.i. From theorem 3.1 the solutions of the system are infinite.

\subsection{Example 2}
Consider the matrix differential equation (5) and let 
\[
F=\left[\begin{array}{ccccc} 1&1&1&1&1\\0&1&1&0&1\\1&1&1&1&1\\0&1&1&0&1\\1&0&1&0&0\\0&0&1&1&1\end{array}\right],
\]
and
\[
G=\left[\begin{array}{ccccc} 1&2&2&1&2\\0&2&2&0&2\\1&2&2&2&3\\0&2&3&1&3\\0&0&0&0&0\\1&0&1&0&0\end{array}\right]
\]
The matrices F, G are non square. Thus the matrix pencil sF-G is singular with invariants, s-2, s-1 the finite elementary divisors, the infinite elementary divisors are of degree 1 and 
$\zeta_1$=0, $\zeta_2$=1 are the r.m.i. From theorem 3.1 there exist non singular matrices
\[
P=\left[\begin{array}{cccccc}1&-1&0&0&0&0\\0&1&0&0&0&0\\-1&0&1&0&0&0\\0&0&0&0&1&0\\0&0&0&0&0&1\\0&-1&0&1&0&0\end{array}\right]
\]
and
\[
Q=\left[\begin{array}{ccccc}0&0&1&1&-1\\1&1&-1&-1&0\\0&0&-1&0&1\\1&0&-1&-1&1\\-1&0&2&1&-1\end{array}\right]
\]
such that  
\begin{equation}
 PFQ=F_K = \left[\begin{array}{ccccc}1&0&0&0&0\\0&1&0&0&0\\0&0&0&0&0\\0&0&0&1&0\\0&0&0&0&1\\0&0&0&0&0\end{array}\right]
\]
\[
PGQ=G_K=\left[\begin{array}{ccccc}1&0&0&0&0\\0&2&0&0&0\\0&0&1&0&0\\0&0&0&0&0\\0&0&0&1&0\\0&0&0&0&1\end{array}\right]
\end{equation}
with
\[
Q_p=\left[\begin{array}{cc}0&0\\1&1\\0&0\\1&0\\-1&0\end{array}\right]
\]
Let the initial values of the system be 
\[
Y(0)=\left[\begin{array}{c}0\\-1\\0\\1\\-1\end{array}\right]
\]
Then
\[
Y(0)\in colspanQ_p
\]
and from theorem 3.1 the solution of the system is
\begin{equation}
Y(t)=\left[\begin{array}{cc}0&0\\1&1\\0&0\\1&0\\-1&0\end{array}\right]\left[
     \begin{array}{cc} e^t&0\\0&e^{2t}\end{array}\right]Z_p(0)
\end{equation}
and by calculating $Z_p(0)$ we get
\[
Y(0)=\left[\begin{array}{cc}0&0\\1&1\\0&0\\1&0\\-1&0\end{array}\right]Z_p(0)
\]
or 
\[
Z_p(0)=
\left[\begin{array}{c}
 1\\-2
\end{array}\right]
\]
and the solution of the system is 
\[
Y(t)=\left[\begin{array}{c} 0\\e^t-2e^{2t}\\0\\e^t\\-e^t\end{array}\right]
\]
Next assume the initial conditions 
\[
Y_0=
\left[\begin{array}{c}  0\\0\\0\\1\\1
\end{array}\right]
\]
Then 
\[
Y_0\notin colspanQ_p
\]
the initial conditions are non consistent and the solution for every $k \geq 0$ is 
\[
Y(t)=\left[\begin{array}{cc} 0&0\\e^t&e^{2t}\\0&0\\e^t&0\\-e^t&0\end{array}\right]C
\]
where C=$\left[\begin{array}{c} c_1\\c_2\end{array}\right]$ is constant and the dimension of the solution vector space of the system is 2.

\section*{Conclusions}

In this article, we study an initial value problem of a class of linear rectangular matrix
differential equations. First by taking into consideration that the
relevant pencil is singular, we decompose the autonomous linear differential system into five
sub-systems and we provide necessary and sufficient conditions for existence and uniqueness of
solutions. As a further extension of the present paper, we can discuss the non-homogeneous case and the case of a singular discrete time system with a singular matrix pencil. For all this there is some research in progress.

\end{document}